\documentclass[12pt]{article}
\usepackage{amstex, amsthm}
\input amssym.def
\input amssym
\title{{\bf Rank one lattice type vertex operator algebras and their
automorphism groups, II: E-series}
\footnotetext{The first author is supported by NSF grant DMS-9700923 and a
research grant from the Committee on Research, UC Santa Cruz.}\footnotetext{
The second author is supported by NSF grant DMS-9623038 and the
University of Michigan Faculty Recognition 
Grant (1993--96).}}
\author{Chongying Dong,  Robert L.~Griess Jr., Alex Ryba}
\date{}

\makeatletter

\begin{document}

\newtheorem{thm}{Theorem}[section]
\newtheorem{prop}[thm]{Proposition}
\newtheorem{lem}[thm]{Lemma}
\newtheorem{rem}[thm]{Remark}
\newtheorem{cor}[thm]{Corollary}
\newtheorem{conj}[thm]{Conjecture}
\newtheorem{de}[thm]{Definition}
\newtheorem{notation}[thm]{Notation}
\pagestyle{plain}
\maketitle

\def\vv{ {\hbox {\bf 1} }}  
\def \Z{\Bbb Z}
\def \F{\Bbb F}
\def \C{\Bbb C}
\def \R{\Bbb R}
\def \Q{\Bbb Q}
\def \N{\Bbb N}
\def \D{{\cal D}}
\def \wt{{\rm wt}}
\def \tr{{\rm tr}}
\def \sp{{\rm span}}
\def \Res{{\rm Res}}
\def \Res{{\rm QRes}}
\def \End{{\rm End}}
\def \E{{\rm End}}
\def \Ind {{\rm Ind}}
\def \Irr {{\rm Irr}}
\def \Aut{{\rm Aut}}
\def \Hom{{\rm Hom}}
\def \mod{{\rm mod}}
\def \ann{{\rm Ann}}
\def \<{\langle} 
\def \>{\rangle} 
\def \t{\tau }
\def \a{\alpha }
\def \e{\epsilon }
\def \l{\lambda }
\def \L{\Lambda }
\def \g{\frak g}
\def \b{\beta }
\def \om{\omega }
\def \o{\omega }
\def \k{\kappa}
\def \c{\chi}
\def \ch{\chi}
\def \cg{\chi_g}
\def \ag{\alpha_g}
\def \ah{\alpha_h}
\def \ph{\psi_h}
\def \pf{{\bf Proof. }}
\def \voa{{vertex operator algebra\ }}
\def \svoa{{super vertex operator algebra\ }}
\def \lc{L_C}
\def \tlc{\widetilde{L}_C}
\def \tv{\widetilde{V}_L}
\def \vlc{V_{L_C}}
\def\tvlc{\widetilde{V}_{L_C}}
\def\vtlc{V_{\widetilde{L}_C}}
\def\tvtlc{\widetilde{V}_{\widetilde{L}_C}}
\def\ha{\frac{1}{2}}
\def\se{\frac{1}{16}}

\def\Ve{V^{0}}
\def\M{\rm \ {\sl  M}\llap{{\sl I\kern.80em}}\ }
\def\xx{\em}
\def\s{\sigma}
\def\a{\alpha} 
\def\t{\tau}
\def\b{\beta} 

\def\xap{x_\a^+}
\def\xam{x_\a^-}

\def\la{\langle}
\def\ra{\rangle}
\def\rtar{\rightarrow}
\def\mt{\mapsto}


{\bf Abstract.} Let $L$ be the $A_1$  root lattice and $G$ a finite subgroup of
$Aut(V)$, where $V=V_L$ is the associated lattice VOA (in this case,
$Aut(V) \cong PSL(2,\Bbb C)$).  The fixed point subVOA, $V^G$ was
studied in [DG], which finds a set of generators and determines the
automorphism group when  $G$ is cyclic (from the ``$A$-series") or
dihedral (from the ``$D$-series").  In the present article, we obtain
analogous results for the remaining possibilities for $G$, that it
belong to the ``$E$-series": $G\cong Alt_4, Alt_5$ or $Sym_4.$

\section{Introduction}

The paper is sequel to [DG] in which we determined a set of generators
and the full automorphism groups of $V_{L_{2n}}^+$ where
$L_{2n}$ is a rank 1 lattice spanned by an element of squared  length $2n$ 
and  $V_{L_{2n}}^+$ is the fixed points of the lattice VOA 
$V_{L_{2n}}$ under an  automorphism of $V_{L_{2n}}$ lifting the $-1$ 
isometry of $L_{2n}.$  
In this paper we determine a set of generators and the full
automorphism group of $V_L^G$ when $L=L_2$ is the root lattice of type 
$A_1$ and
$G$ is an automorphism group of type tetrahedral, octahedral or icosahedral.

The graded dimensions of $V_{L_{2n}}$,  $V^+_{L_{2n}}$, and the three
$V^G_{L_{2n}}$ realize all the partition functions of rank 1 rational
conformal field theories; such functions (but not the VOAs themselves)
are classified in physics literature [G] and [K].  It is unknown
whether two inequivalent rank 1 rational VOAs may have the same graded
dimension and it is also unknown whether all the VOAs above are rational. 
Certainly, the $V_{L_{2n}}$ are rational (see [D] and [DLM2])
 and some progress has
been made towards showing that the $V^+_{L_{2n}}$ are rational,  namely
the finiteness of the number of isomorphism types of irreducible
modules has been satisfied [DN1]. 

It is well known that
the finite subgroups of $SO(3)$ are labeled by
the simply-laced Lie algebras. If $G$ is of type $A$ or $D,$
$V_{L_2}^G$ is $V_{L_{2n}}$ or $V_{L_{2n}}^+$ for 
some $n.$ Since the full automorphism
groups for all lattice vertex operator algebras and $V_{L_{2n}}^+$ have
been determined in [DN2] and [DG],  the results in this paper 
complete the determination of generators and full automorphism
groups for this set of  vertex operator algebras 
of rank 1. Using the results from [DN2], [DG] and the present paper
one can easily see that the set of isomorphism types 
$${\cal S}=\{V_{L_{2n}},V_{L_{2n}}^+, V_{L_2}^G
|n\geq 1, G=Alt_4,Sym_4,Alt_5\}$$
is closed in the sense that for any $V\in {\cal S}$ and a finite 
subgroup $G$ of
$Aut(V),$ then $V^G\in {\cal S}.$ 

The paper is organized as follows. In Section 2, we review the invariant 
theory for the subgroups of $SO(3)$ of $E$-series following [S]. In Section 3
we determine the generators and the automorphism groups of $V_{L_2}^G$ 
for $G\cong Alt_4, Sym_4,Alt_5$ which are the subgroups of $SO(3)$ of type $E.$

We assume that the reader has some familiarity  with the definition
of vertex operator algebra and vertex operator algebras associated to
even positive definite lattices as presented in [B], [FLM] and [DG].

\section{Representations of $SL(2,\C)$}
\setcounter{equation}{0}

\def\dim{{\rm dim}}  
\def\rank{{\rm rank}}  
\def\F#1{{\bf F_{#1}}}  
\def\gg{SL_2(\C)}             
\def\gz{SL_2(\Z)}             
\def\gr{SL_2(R)}              
\def\gs{SL_2(R^*)}            
\def\gq#1{SL_2(\F{#1})}      
\def\Wm#1{W_{#1}}            
\def\La#1{\Lambda_{#1}}      
\def\oWm#1{\overline{W_{#1}}}
\def\proj#1{p_{#1}}                
\def\oproj#1{\overline{p_{#1}}}    
\def\Symd{S}           
\def\Altd{T}           
\def\Alte{I}           

\bigskip

The following notation will be used throughout the paper.
\halign{#\hfil&\quad#\hfil\cr

$\Wm{m}$ & The $m$-dimensional irreducible module for $\gg$. \cr

$\proj{m}$ & The projection from a finite dimensional $\gg$-module
onto its \cr

& $\Wm{m}$-homogeneous component\cr

$\Altd$ & A particular copy of the finite group $2Alt_4$ in $\gg$.\cr

$\Symd$ & A particular copy of the finite group $2Sym_4$ in $\gg$.\cr

$\Alte$ & A particular copy of the finite group $2Alt_5$ in $\gg$.\cr

$V^H$ & The fixed points of the action of a group $H$ on a module $V$.\cr
}
\bigskip

We recall that the group $\gg$ has a unique irreducible module
$\Wm{m}$ of any finite dimension $m$ [H].  This module contains an
$m$-dimensional integral representation, spanned by a Chevalley basis,
of the integral Chevalley group $\gz.$ We shall  
write this integral representation as $\La{m}$.  
We shall need to work with $SL_2(A)$ modules for various
choices of a subring $A$ in $\C$.  In particular, we write
$\Wm{m,A}$ for the $SL_2(A)$-module $A\otimes_{\Z}\La{m}$.
We write $R$ for a ring of
algebraic integers, and we let bars indicate reduction
modulo a prime containing $p$ in $R$ and for the
result of tensoring with $\bar R.$ 

\smallskip

We shall be interested in tensor products of pairs of $\gg$-modules.  The
decompositions of these tensor products into irreducibles are given by the
Clebsch-Gordan formula:
$\Wm{m} \otimes \Wm{n} = \Wm{m+n-1} \oplus \Wm{m+n-3} \oplus \ldots \oplus \Wm{n - m +1}$,
which holds whenever $n \geq m$ [H]. 
A similar decomposition of $KSL_2(K)$ modules
holds for any subfield $K$ in $\C$. 
The main result of this section is Theorem
\ref{t2.1}, which is needed in Section 3.

\begin{thm}\label{t2.1} We have

\noindent(i) $\proj{19} (\Wm9^\Symd \otimes \Wm{13}^\Symd) \neq 0$. 

\noindent(ii) $\proj{31} (\Wm{13}^\Alte \otimes \Wm{21}^\Alte) \neq 0$. 

\noindent(iii) The 1-dimensional spaces $\proj{13} (\Wm{7}^\Altd \otimes \Wm{7}^\Altd)$ and
               $\proj{13} (\Wm{7}^\Altd \otimes \Wm{9}^\Altd)$ are distinct.

\noindent(iv) The 1-dimensional spaces $\proj{13} (\Wm{9}^\Altd \otimes \Wm{9}^\Altd)$ and
               $\proj{13} (\Wm{7}^\Altd \otimes \Wm{9}^\Altd)$ are distinct.

\end{thm}

To establish this theorem, we must establish lower bounds for projections
of particular subspaces of tensor products of $\gg$-modules.   (The upper bounds
of dimension 1
implicit in (iii) and (iv) are immediate consequences of the Clebsch-Gordan
formula.)  We shall establish
these claims by performing explicit computations in analogous modules for
a finite group of type $SL_2(\bar R)$ and lifting the results to characteristic
0.  We begin by obtaining conditions under which we can
lift statements about the dimension of images of projection maps.

\smallskip

\begin{lem}
Suppose 
that $L$ and $M$ are finite rank $R$-torsion free 
$RSL_2(R)$-modules
that are equivalent over the field of fractions of $R$.
Then $\bar L$ and $\bar M$ are $\bar R SL_2(\bar R)$-modules
with identical sets of composition factors.
\end{lem}

Proof.
The first argument for 82.1 in [CR] (which
establishes an analogous result for representations of a finite
group) applies without change.
\qed

\smallskip

Let $\Gamma_{m,n}$ denote the set of degrees of irreducible
constituents in the Clebsch-Gordan decomposition.  Thus,
$\Gamma_{m,n} = 
\lbrace m + n - 1, m + n - 3, \ldots, m + n + 1 - 2 \ min(m,n) \rbrace$.
When $k \in \Gamma_{m,n}$ we write $q_k$ for the composite of 
the map given
by extending the scalars for $\Wm{m,R} \otimes \Wm{n,R}$ to the
field of fractions, $K$ of $R$, followed by the 
projection 
onto the $k$-dimensional irreducible constituent of the
$KSL_2(R)$-module.  

\smallskip

\begin{lem}
Suppose that $k \in \Gamma_{m,n}$,
that $\overline{\Wm{m,R} \otimes \Wm{n,R}}$ is
completely reducible and that the degrees of its irreducible
constituents are given by $\Gamma_{m,n}$. 
Then $\overline{Im(q_k)}$ is irreducible and
$\overline{q_k} : \overline{\Wm{m,R} \otimes \Wm{n,R}}
\rightarrow \overline{Im(q_k)}$ is a surjection.
\end{lem}

Proof.
Let $K$ be the field of fractions of $R$.
The $R$-free $RSL_2(R)$-modules
$\Wm{m,R} \otimes \Wm{n,R}$ and
$\bigoplus_{k \in \Gamma_{m,n}} \Wm{k,R}$
both extend to $KSL_2(R)$-modules isomorphic to
$\Wm{m,K} \otimes \Wm{n,K}$.  Thus, according to
Lemma 2.2, the modular reductions
$\overline{\Wm{m,R} \otimes \Wm{n,R}}$ and
$\bigoplus_{k \in \Gamma_{m,n}} \overline{\Wm{k,R}}$
have the same sets of composition factors.  However,
our hypothesis about the composition
factors of the first of these modules shows that
$\overline{\Wm{k,R}}$ is an irreducible $\bar{R}SL_2(\bar{R})$-module
of degree $k$.  

Now, ${Im(q_k)}$ is an $RSL_2(R)$-module
which is $K$-equivalent to $\Wm{k,R}$.  Hence, by
Lemma 2.2, 
$\overline{Im(q_k)}$ is also an irreducible 
$\bar{R}SL_2(\bar{R})$-module of degree $k$.  
Since $q_k$ is a surjection from
$\Wm{m,R} \otimes \Wm{n,R}$ to $Im(q_k)$, 
$\overline{q_k}$, the result of tensoring with $\bar R$,
is also a surjection.
\qed

\bigskip

\begin{cor}
Suppose that $k \in \Gamma_{m,n}$, and
that $\overline{\Wm{m,R}} \otimes \overline{\Wm{n,R}}$ is
completely reducible and that the degrees of its irreducible
constituents are given by $\Gamma_{m,n}$. 
Suppose that $S =
\lbrace (m_1,n_1), (m_2, n_2), \ldots, (m_r,n_r)\rbrace
\subset \Wm{m,R} \times \Wm{n,R}$ has the property
that 
$\overline{m_1}\otimes\overline{n_1}, \overline{m_2}\otimes\overline{n_2}, 
\ldots, \overline{m_r}\otimes\overline{n_r}$ have linearly independent
images under an $\bar{R}SL_2(\bar{R})$-module
homomorphism from
$\overline{\Wm{m,R}} \otimes \overline{\Wm{n,R}}$ onto its
irreducible image of degree $k$.  Then 
$p_k(m_1\otimes n_1), p_k(m_2\otimes n_2), 
\ldots, p_k(m_r\otimes n_r)$ are linearly independent.
\end{cor}

Proof.
The modules
$\overline{\Wm{m,R}} \otimes_{\bar R} \overline{\Wm{n,R}}$ and
$\overline{\Wm{m,R} \otimes_R  \Wm{n,R}}$ are naturally isomorphic.
Thus, we can identify the essentially unique
$\bar{R}SL_2(\bar{R})$-module
homomorphism from
$\overline{\Wm{m,R}} \otimes \overline{\Wm{n,R}}$ onto its
irreducible image of degree $k$ with the map
$\overline{q_k}$ of Lemma 2.3.  Independence of images
under $\overline{q_k}$ implies independence of the
corresponding images under $q_k$ and $p_k$.
\qed

\bigskip

Suppose that $F$ is a finite subgroup of $\gg$ and that the
corresponding character of $F$ can be written over a ring
of integers 
$R$ with $\bar R = \F{p}.$ 
Moreover, if the tensor product $\oWm{n} \otimes \oWm{m}$
of $SL_2(\bar R)$-modules meets the conditions of Corollary 2.4,
then the $p$-modular reduction of
$\proj{k}(\Wm{n}^F \otimes \Wm{m}^F)$ is the projection
of $\oWm{n}^F \otimes \oWm{m}^F$ onto its unique $k$-dimensional
summand.  Our strategy for proving Theorem \ref{t2.1} is to 
compute the latter projections.  Such work with
barred objects is relatively pleasant since it 
involves linear algebra over the integers modulo a prime.  

We write $r_k$ for the projection of an
$SL_2(\bar R)$-module onto an irreducible summand of degree $k$ when such a
summand exists and that irreducible has multiplicity $1$ in the module
(the latter condition implies uniqueness of such a projection).

Here is the computation (which says in
part (ii) that ``$r_k=\overline{p_k},$'' for suitable $k$).
Theorem \ref{t2.1} follows.

\begin{prop} Let $p = 101$, then: 

\noindent(i)  The traces for elements  of
$\Altd$, $\Symd$, and $\Alte$ in $\gg$ reduce in $\bar R$ 
 to elements
of the prime field $\F{101}$. So if $Z=T,S$ or $I,$  we may
assume that $R$ satisfies   $Z\leq SL_2(R)$
and  $\bar R={\bf F}_{101}.$

\noindent(ii) The modules 
$\oWm{9} \otimes \oWm{13}$, 
$\oWm{13} \otimes \oWm{21}$, 
$\oWm{7} \otimes \oWm{7}$, 
$\oWm{7} \otimes \oWm{9}$, and 
$\oWm{9} \otimes \oWm{9}$ for $SL(2,101)$ are all completely reducible
and the degree sets of their irreducible summands are
$\lbrace 5, 7, 9, 11, 13, 15, 17, 19, 21\rbrace$,
$\lbrace 9, 11, 13, 15, 17, 19, 21, 23, 25, 27, 29, 31, 33\rbrace$,
$\lbrace 1, 3, 5, 7, 9, 11, 13\rbrace$,
$\lbrace 3, 5, 7, 9, 11, 13, 15\rbrace$,
and $\lbrace 1, 3, 5, 7, 9, 11, 13, 15, 17\rbrace$.  These degree sets
match the degrees in the corresponding versions of the
Clebsch-Gordan formula.

\noindent(iii) $r_{19}
(\oWm9^\Symd \otimes \oWm{13}^\Symd) \neq 0$. 

\noindent(iv) $r_{31}
(\oWm{13}^\Alte \otimes \oWm{21}^\Alte) \neq 0$. 

\noindent(v) The 1-dimensional spaces $r_{13}
(\oWm{7}^\Altd \otimes \oWm{7}^\Altd)$ and
$r_{13}(\oWm{7}^\Altd \otimes \oWm{9}^\Altd)$ are distinct.

\noindent(vi) The 1-dimensional spaces $r_{13}
(\oWm{9}^\Altd \otimes \oWm{9}^\Altd)$ and
$r_{13}(\oWm{7}^\Altd \otimes \oWm{9}^\Altd)$ are distinct.

\end{prop}

{\bf Method.} For (i), we note that the only character irrationalities that we
need to 
consider involve $5^{th}$ roots of unity, and that ${\bf F}_{101}$ contains
fifth roots of unity. Therefore, we may work over a suitable ring of integers, $R,$ such that $\bar R={\bf F}_{101}.$

For (ii), we compute the decompositions of the tensor products with
the Meat-Axe [Pa].  Moreover, by using the Meat-Axe to determine bases
for the irreducible submodules of the dual spaces of these tensor
products, we obtain matrices representing all projection maps,
$r_k$, onto irreducible summands of the tensor products.

For (iii), (iv), (v) and (vi), we compute fixed point spaces under
finite subgroups of $SL_2(\bar R)$ with the Meat-Axe.  The projections of
images of tensor products of these spaces are then determined by
using the representations of the projection maps that we computed in (ii). 
\qed 

\bigskip

\section{Generators and automorphisms of $V_{L_2}^G$}

\def\RR{{\Bbb R}}
\def\CC{{\Bbb C}}
\def\NN{{\Bbb N}}
\def\ZZ{{\Bbb Z}}
\def\FF{{\Bbb F}}
\def\QQ{{\Bbb Q}}
\def\AA{{\Bbb A}}
\def\BB{{\Bbb B}}
\def\DD{{\Bbb D}}
\def\EE{{\Bbb E}}
\def\GG{{\Bbb G}}
\def\HH{{\Bbb H}}
\def\II{{\Bbb I}}
\def\JJ{{\Bbb J}}
\def\KK{{\Bbb K}}
\def\LL{{\Bbb L}}
\def\MM{{\Bbb M}}
\def\OO{{\Bbb O}}
\def\PP{{\Bbb P}}
\def\SS{{\Bbb S}}
\def\TT{{\Bbb T}}
\def\UU{{\Bbb U}}
\def\VV{{\Bbb V}}
\def\WW{{\Bbb W}}
\def\XX{{\Bbb X}}
\def\YY{{\Bbb Y}}



\def\App#1#2{{\bf (#1.#2) Application. }}
\def\Hyp#1#2{{\bf (#1.#2)  Hypothesis. }}
\def\Not#1#2{{\bf (#1.#2)  Notation. }} 
\def\Def#1#2{{\bf (#1.#2) Definition. }}
\def\Lem#1#2{{\bf (#1.#2)  Lemma. }}
\def\Th#1#2{{\bf (#1.#2)  Theorem. }}
\def\Prop#1#2{{\bf (#1.#2)  Proposition. }}
\def\Rk#1#2{{\bf (#1.#2) Remark. }}
\def\Ex#1#2{{\bf (#1.#2) Example. }}


\def\Euc{$\RR^{\ell}$ }
\def\Pf{{\bf Proof. }}
\def\pa{\vskip .6cm} 

\def\Eh{$E_8(\CC)$}
\def\Eg{$E_7(\CC)$} 
\def\tEg{$2E_7(\CC)$}
\def\Ef{$E_6(\CC)$}
\def\tEf{$3E_6(\CC)$}

\def\Ehh{$E_8(\CC)~$}
\def\Gb{$G_2(\CC)$}
\def\Gbb{$G_2(\CC)~$}
\def\Fd{$F_4(\CC)$}
\def\Fdd{$F_4(\CC)~$}

\def\g{$\frak g$} 
\def\gg{$\frak g$ }

\def\Vnat{$V^\natural$}
\def\Vnatt{$V^\natural$ }

\def\half{$1 \over 2$}

\def\tCo1{$2^{24}.Co_1$}  
\def\ctCo1{$2^{1+24}_+ \cdot Co_0$}

\def\tc1{$2^{24}\cdot Co_1$}  
\def\ctc1{$2^{1+24}_+ \cdot Co_0$}

\def\c1{$Co_1$}
\def\c0{$Co_0$}

\def\L{\Lambda}

\def\a{\alpha}
\def\b{\beta}
\def\g{\gamma}
\def\s{\sigma}
\def\t{\tau}
\def\o{\omega}
\def\O{\Omega}

\def\kd#1#2{\delta _{#1,#2}}

Let $L_2:=\ZZ \a$,  $(\a, \a)= 2$, the $A_1$-lattice, $V:=V_{L_2}$, the
VOA of lattice type based on $L_2$.  
We let $G$ be a subgroup of $Aut(V) \cong PSL(2,\CC)$ isomorphic to $Alt_4, 
Sym_4$ or $Alt_5$.  The irreducible $W_m$ for $SL(2,\CC)$ of
dimension $m$ may be interpreted as a module for $PSL(2,\CC)$ when $m$ is odd.

\begin{thm}\label{t3.1} $Aut(V^G)$ is the identity if
$G\cong Sym_4$ or $Alt_5$ and is isomorphic to $\Z_2$ if $G\cong Alt_4.$ So, in all cases,
$Aut(V^G)  \cong N_{Aut(V)}(G)/G$
\end{thm}

It is well-known that $SL(2,\CC)$ acts on $\CC[x,y]$ as algebra isomorphism
such that $\CC x+\CC y$ is a natural $SL(2,\CC)$-module. One can identify 
$W_{m+1}$ with the space of degree $m$ homogeneous polynomials. As an 
$SL(2,\CC)$-module, $\C[x,y]$ has a decomposition
$$\C[x,y]=\bigoplus_{m\geq 1}W_m.$$
Let $\tilde G$ be the preimage of $G$ in $SL(2,\CC)$ (we may assume
$\tilde G=T,S$ or $I$) and 
$A$ the algebra of invariants for the action of $\tilde G$ on $\CC [x,y].$
The following proposition can be found in [S]. In all cases, $A$ is a 
quotient of a 
polynomial ring with $3$ generators modulo an ideal generated by 
a single relator, indicated below. 

\begin{prop}\label{p3.2} The algebra has a set of generators as
follows  ($f_n, g_n$ etc. denote homogeneous polynomials of degree $n$):

\noindent(i)   If $\tilde G \cong SL(2,3),$ $A=\CC[f_6,f_8,f_{12}],$ 
subject to the relation $f_6^4+f_8^3+f_{12}^2=0.$

\noindent(ii) If $\tilde G \cong 2\cdot Sym_4,$ $A=
\CC[ g_{8}, g_{12}, g_{18} ],$ subject to the relation $g_{18}^2+g_8^3g_{12}+g_{12}^3=0.$

\noindent(iii) If $\tilde G \cong SL(2,5),$ $A=
\CC[ h_{12},h_{20},h_{30}],$ subject to the relation 
$h_{12}^5+h_{20}^3+h_{30}^2=0.$
\end{prop}

Let $L(c,h)$ be the irreducible highest weight module for the Virasoro
algebra with central charge $c$ and highest weight $h$ for $c,h\in \C.$ 
The subspace $H$ of highest weight vectors for $Vir$ in $V$  
is isomorphic to the subspace $\CC [x,y]^+$ of even polynomials in
$\CC [x,y]$, and  we may and do assume that this isomorphism
$\varphi : H \rightarrow \CC [x,y]^+$ 
preserves degree.  
We identify $W_{2m+1}$ with the degree $m$ part of $H$.  

Recall from [DG] that there is an
isomorphism  $$ V \cong \sum_{m \ge 0} W_{2m+1} \otimes
L(1,m^2), \leqno (*)$$
as modules for $SL(2,\C) \times L(1,0).$  Let $\pi _m$ be the
projection of $V$ to $W_{2m+1} \otimes
L(1,m^2)$.  Note that since $dim(W_1)=1$, $\pi_1$ can be interpreted as a
map onto $L(1,0)$.  

We need the following result from [DM2] (see Lemma 2.3; also see [DM1]):
\begin{lem}\label{l3.3} Let $K$ be a compact Lie group which acts 
continuously on a vertex operator algebra $U.$ Let $M,N$ be two
finite dimensional $K$-submodule of $U.$ Then there exists $n$ such that
the linear span of $\sum_{m\geq n}s_mt,$ for $s\in M,t\in N,$ is
isomorphic to $M\otimes N$ as $K$-modules. 
\end{lem}

For convenience we set $v^1= \varphi^{-1} (f_6)$ (resp. $g_8$ or $h_{12}$),
$v^2=  \varphi^{-1} (f_8)$ (resp. $g_{12}$ or $h_{20}$) and
$v^3=  \varphi^{-1} (f_{12})$ (resp. $g_{18}$ or $h_{30}$) 
if $G\cong Alt_4$ (resp. $Sym_4$ or
$Alt_5$).  Then we have 
\begin{prop}\label{p3.4} 
The vertex operator algebra $V^G$ is generated by $\{ \o,  v^1,v^2,v^3 \}$.
\end{prop}

\pf We prove the result 
for $\tilde G=T$ the other cases being similar.
 First note that the algebra $\C[x,y]^{T}$ is generated by 
$f_6, f_8$ and $f_{12}$ and the algebra product can factor through 
the tensor product. The $T$-invariants of $V$ have the form 
$$V^{T}=\bigoplus_{m\geq 0}W_{2m+1}^{T}\otimes L(1,m^2).$$
It is enough to show that if $W_{2s+1}^{T}\otimes L(1,s^2)$
and $W_{2t+1}^{T}\otimes L(1,t^2)$ can be generated by 
$\o$ and the $v^i$ then so are $W_{2(s+t)+1}\otimes L(1,(s+t)^2)$.
We assume that $s\geq t.$

By Lemma \ref{l3.3}, span$\{u_mv|u\in W_{2s+1}\otimes L(1,s^2),
v\in W_{2t+1}\otimes L(1,t^2), m\in \Z\}$ is exactly
the subspace
$$\bigoplus_{l=2s-2t+1, 2s-2t+3,...2s+2t+1}W_{l}\otimes L(1,(\frac{l-1}{2})^2)$$
and span$\{u_mv|u\in W_{2s+1}^{T}\otimes L(1,s^2),
v\in W_{2t+1}^{T}\otimes L(1,t^2), m\in \Z\}$ is exactly
the subspace
$$\bigoplus_{l=2s-2t+1, 2s-2t+3,...2s+2t+1}p_l(W_{2s+1}^{T}\otimes W_{2t+1}^{T})\otimes L(1,(\frac{l-1}{2})^2).$$ 
Since $p_{2(s+t)+1}(W_{2s+1}^{T}\otimes W_{2t+1}^{T})=W_{2(s+t)+1}^{T},$ 
we immediately see that $W_{2(s+t)+1}\otimes L(1,(s+t)^2)$ can be generated
by $W_{2s+1}^{T}\otimes L(1,s^2)$
and $W_{2t+1}^{T}\otimes L(1,t^2).$ As a result, $W_{2(s+t)+1}\otimes L(1,(s+t)^2)$ can be generated by $\o$ and $v^i.$ 
\qed 

\begin{lem}\label{l3.3'} For any $k\geq 1,$ there is an invariant 
 bilinear form on $W_{k}$, and for any $x, y \in W_{k}$ which are not orthogonal, $p_1(x \otimes y) \ne 0.$
This applies to $x=y\ne 0$ in  
$W_{k}^G$ whenever $dim(W_{k}^G)=1.$
\end{lem}

\pf Since $W_k$ is the unique irreducible of dimension $k,$ it is self-dual 
module. A nonzero invariant bilinear form is unique up to
scalar. Since $p_1$ maps $W_k\otimes W_k$ onto $W_1,$ this projection
is essentially that bilinear form. Since the subspace $W_k^G$ is nonsingularly
paired with itself under such a bilinear form, the last statement
follows. \qed

{\bf Proof of Theorem \ref{t3.1}:} (i) Let $G\cong Sym_4.$  
Let $\s \in Aut(V^G)$. 
Since $\s$ preserves weights and fixes
$\o$, it stabilizes each $W_{2m+1}^{G}$ which is the subspace of
highest weight vectors of weight $m^2$ in $V^{G}$ for the Virasoro 
algebra. By Proposition \ref{p3.2}, 
$W_{9}^{G}=\C g_8, W_{13}^{G}=\C g_{12}$ and $W_{19}^{G}=\C g_{18}.$
Thus, there are scalars $c_i \in \C$ such that  
$\s g_i = c_i g_i$ for $i=8,12,18.$ By Lemmas \ref{l3.3} and \ref{l3.3'}
there exists $m_i\in \Z$ such that $0\ne \pi_1((g_i)_{m_i}g_i)\in L(1,0).$
Since $\s$ is trivial on $L(1,0),$ we immediately have
$c_i^2=1$ for all $i.$ That is, $c_i=\pm 1.$  

{}From Proposition \ref{p3.2}, $W_{37}^{G}$ is $2$-dimensional. Since
$g_{18}^2\in W_{37}^{G}$ which is regarded as a subspace of $\C[x,y],$  
we see from
Lemma \ref{l3.3} there exists $m\in \Z$ such that 
$0\ne \pi_{6}((g_{18})_mg_{18})\in W_{37}^{G}\otimes L(1,18^2).$ So,
$\s$ has an eigenvalue 1 on $W_{37}^{G}\otimes L(1,18^2).$
Using the actions of the Virasoro operators $L(n)$ for $n\geq 0$
on  $\pi_{6}((g_{18})_mg_{18})$ we get an eigenvector in $W_{37}^{G}$
for $\sigma$ with eigenvalue 1. 
Similarly, there exist $n_1,n_2, s_1,s_2,s_3\in \Z$ such that $\pi_{37}((g_{12})_{n_1}
\pi_{25}((g_{12})_{n_2}g_{12}))\in W_{37}^{G}\otimes L(1,18^2)$
 is an eigenvector
of $\s$ with eigenvalue $c_{12}^3$ and
$\pi_{37}((g_{8})_{s_1}\pi_{27}[(g_{8})_{s_2}\pi_{21}[(g_{8})_{s_3}g_{12}]])$
is an eigenvector of $\s$ with eigenvalue $c_8^3c_{12}.$  As a result 
$W_{37}^{G}$ contains 3 eigenvectors $u,v,w$ of $\s$ with eigenvalues
$1, c_{12}^3$ and $c_8^3c_{12}.$ 
From Proposition
\ref{p3.2}, any two of $\{u,v,w\}$ are linearly
independent and the three are linearly dependent. 
Since $W_{37}^{G}$ is $2$-dimensional, we conclude
$c_{12}^3=c_8^3c_{12}=1,$  whence
$c_8=c_{12}=1.$

It remains to show that $c_{18}=1.$ By Theorem \ref{t2.1} (i)
and Lemmas \ref{l3.3},  there exists $m\in \Z$ such that 
$\pi_{19}((g_{8})_mg_{12})\in W_{19}^{G}\otimes L(1,81)$ is an eigenvector
of $\s$ with eigenvalue $1.$  Since $W_{19}^{G}$ is 1-dimensional
we see that $\s g_{18}=g_{19}.$ Since
$V^{G}$ is generated by $\o$ and $g_i$ (see Proposition \ref{p3.4}),
we conclude that $\s$ is the identity map.

(ii) The proof in the case that $G\cong Alt_5$ is similar to that in case (i).

(iii) Let $G\cong Alt_4.$  In this case $V^{G}$ is generated by $\o$ and
$f_6,f_8,f_{12}$ by Proposition \ref{p3.4}. Let $\s\in Aut(V^{G}).$
Since $W_{7}^{G}=\C f_6$ and $W_{9}^{G}=\C f_8$ there exist
$c_6,c_8\in\C$ so that $\s f_6=c_6 f_6$ and $\s f_8=c_8 f_8.$
Using the same argument as used in the proof of case (i) gives
$c_6=\pm 1$ and $c_8=\pm 1.$ Note that $W_{13}^{G}=\C f_6^2+\C f_{12}$
is 2-dimensional.  From Theorem \ref{t2.1} (iii) and Lemma
\ref{l3.3}, there exist $m,n\in \Z$ such that 
$\pi_{13}((f_6)_mf_6),\pi_{13}((f_6)_nf_8)\in W_{13}^{G}\otimes L(1,36)$
are linearly independent eigenvectors of $\s$ with eigenvalues
$1$ and $c_6c_8.$ This implies that $V^{G}$ is generated by $\o,$
$f_{6}$ and $f_8$ and the automorphism group is isomorphic to one
of $1,$ $\Z_2$ and $\Z_2\times \Z_2.$  

Recall that the normalizer $N(G)=N_{Aut(V)}(G)$ is isomorphic to $Sym_4.$
{}From Lemma 3.2 of [DM1], we know that $V^{N(G)}$ and $V^{G}$ are
different, whence an element of $N(G)/G \cong \Z_2$ 
acts on $V^{G}$ as a nontrivial
automorphism, denoted by $\s.$ So, $Aut (V^{G})$ is either $\Z_2$ or
$\Z_2\times \Z_2.$ If $Aut(V^{G})$ is isomorphic to $\Z_2\times \Z_2,$
then there exists $\tau\in Aut (V^{G})$ such that $Aut(V^{G})$ is
generated by $\s$ and $\tau.$ By Theorem 2.4 of [DLM1], $V^{G}$ can be
decomposed
$$V^{G}=V^{G}_{1,1}\bigoplus V^{G}_{1,-1}\bigoplus V^{G}_{-1,1}\bigoplus V^{G}_{-1,-1}$$
where $ V^{G}_{\mu,\l}=\{v\in V^{G}|\s v=\mu,\tau v=\l v\}$ and
each is nonzero. Moreover $Y(u,z)v\in V^{G}_{\mu_1\mu_2,\l_1\l_2}[[z,z^{-1}]]$ if $u\in V^{G}_{\mu_1,\l_1}$ and $v\in V^{G}_{\mu_2,\l_2}.$ 
It is easy to see that $V^{N(G)}=V^{G}_{1,1}\bigoplus V^{G}_{1,-1}$
and $V^{N(G)}$ has a nontrivial automorphism. This is a contradiction 
to (i).
Thus, $Aut (V^{G})$ must be isomorphic to $\Z_2.$ This completes the proof.
\qed

\begin{rem}{\rm From the proof of Theorem \ref{t3.1} we see, in fact, that
$V^{Alt_4}$ is generated by $\o,$ $f_6$ and $f_8.$ This strengthens 
the result in Proposition \ref{p3.4}.}
\end{rem}

\medskip
\noindent Department of Mathematics, University of California, 
Santa Cruz, CA 95064 USA. Email address: dong@cats.ucsc.edu (C.D.)

\noindent Department of Mathematics, University of Michigan, 
Ann Arbor, MI 48109-1109 USA. Email Address: rlg@math.lsa.umich.edu (R.L.G.)

\noindent Department of Mathematics, Marquette University, Milwaukee, WI 53201-1881 USA. Email address: ryba@math.lsa.umich.edu (A.R.)
\end{document}